\documentclass[a4paper,11pt]{amsart}

\usepackage{amsmath}
\usepackage{amsfonts}
\usepackage{amssymb}
\usepackage{amsthm}
\usepackage[english]{babel}
\usepackage[latin1,ansinew]{inputenc}
\usepackage{a4wide}
\usepackage{dsfont}
\usepackage[active]{srcltx}	

\usepackage{graphicx}
\newtheoremstyle{note}
 	{3pt}
 	{3pt}
 	{\itshape}
 	{}
 	{\bfseries}
 	{.}
 	{.5em}
 	{\thmname{#1}\thmnumber{ \textup{#2}}\thmnote{ \upshape#3}}
\theoremstyle{note}
\newtheorem{thm}{Theorem}
\newtheorem{lemma}{Lemma}

\newtheorem*{thm*}{Theorem}
\newtheorem{defi}{Definition}
\newcommand{\R}{{\mathbb R}}

\newcommand{\nn}{\nonumber}
\newcommand{\bea}{\begin{eqnarray}}
\newcommand{\eea}{\end{eqnarray}}
\newcommand{\beann}{\begin{eqnarray*}}
\newcommand{\eeann}{\end{eqnarray*}}
\newcommand{\ba}{\begin{array}}
\newcommand{\ea}{\end{array}}

\newcommand{\beq}{\begin{equation}}
\newcommand{\eeq}{\end{equation}}
\newcommand{\be}{\begin{equation}}
\newcommand{\ee}{\end{equation}}

\providecommand{\norm}[1]{\lVert#1\rVert}

\title{Standing Solitary Euler-Korteweg Waves are Unstable}
\author{Johannes H\"owing}
\email{Johannes.Hoewing@math.uni-hamburg.de}
\address{University of Hamburg, Department of Mathematics, Germany}
\date{January 13, 2013}
\begin{document}

\maketitle

\section{The result}

The Euler-Korteweg system is given by the equations 
\be\begin{aligned}\label{EKL_2}
V_t-U_y&=0,\\
U_t+p(V)_y&=-(\kappa(V)V_{yy}+\frac{1}{2}(\kappa(V))_yV_y)_y,
\end{aligned}
\ee
with $\kappa(V)>0.$ System \eqref{EKL_2}, and notably its solitary waves 
\be\nn
\begin{pmatrix}
V\\
U
\end{pmatrix}(x,t) =
\begin{pmatrix}
v\\
u
\end{pmatrix}(x-ct)
\quad\text{with}\quad
\begin{pmatrix}
v\\
u
\end{pmatrix}(\pm\infty) = \begin{pmatrix}
v_*\\
u_*
\end{pmatrix},
\ee
have been intensely studied by Benzoni-Gavage et al.\
in \cite{B,BDD,BDDJ}.
\begin{defi}\label{definition1_B}\cite{BDDJ}
A traveling wave $(v,u)$ of \eqref{EKL_2} is called orbitally stable if for each 
$\varepsilon>0,$ there exists a $\delta>0$ such that for any solution
$(V,U)\in (v,u)+C([0,T);H^3(\R)\times H^2(\R))$ of \eqref{EKL_2}, 
closeness at initial time, 
$$
\norm{(V,U)(\cdot,0)-(v,u)(\cdot)}_{H^1\times L^2}<\delta
$$
implies closeness at any time 
$$
\inf_{\sigma\in\R}\norm{(V,U)(\cdot,t)-(v,u)(\cdot+\sigma)}_{H^1\times L^2}<\varepsilon
\quad\text{for all }t>0.
\label{foranyt}
$$
\end{defi}
The following is the point of this short note.
\begin{thm}
All non-trivial standing solitary Euler-Korteweg waves are not orbitally stable.
\end{thm}

\section{The proof}
For fixed base state $v_*$, the solitary waves homoclinic to $v_*$ occur in families 
$(u^c,v^c)$ parametrized by their speed $c$. The proof of Theorem 1 is based on 
\begin{lemma} \cite{B}
A solitary wave $(u^{c_*},v^{c_*})$ is orbitally unstable if 
the moment of instability  
\be\nn\begin{aligned}
m(c)&=\int_{-\infty}^{\infty} \kappa(v)v'^2\;d\xi
\end{aligned}
\ee 
is not convex at $c=c_*$.
\end{lemma}

Theorem 1 follows from 

\begin{lemma}
$m''(0)<0.$ 
\end{lemma}
To prove Lemma 2, we recall that with 
\be
F(v,c)=
-f(v)+f(v_*)-p(v_*)(v-v_*)+\frac{1}{2}c^2(v-v_*)^2,
\quad
-\frac{df(v)}{dv}=p(v),
\ee
the profile equation 
\be\nn
\begin{aligned}
\kappa(v)v''+\frac{1}{2}(\kappa(v))'v'
  &= -\frac{\partial F(v,c)}{\partial v} 
\end{aligned}
\ee
possesses (cf.\ \cite{BDDJ}) a first integral given by
\be
I(v,v')
= \frac{1}{2}\kappa(v)v'^2 + F(v,c).
\ee

As 
\be\nn
m(c) = 2\int_{v_*}^{v_m(c)}\kappa(v)v'\;dv
\ee 
with $v_*,v_m(c)>v_*$ consecutive zeros of $F(\cdot,c)$.
Since $I(v,v')\equiv0$ along solutions, we have
\be\nn
\begin{aligned}
m(c)&=2\int_{v_*}^{v_m(c)}\left(\kappa(v)\right)^{1/2}\,\left(
-2F(v,c)\right)^{1/2}\;dv\\
&= 4\int_0^{(v_m(c)-v_*)^{1/2}} \left(\kappa(v_m(c)-w^2)\right)^{1/2}
\left(-2F(v_m(c)-w^2,c) \right)^{1/2}\,w\;dw,
\end{aligned}\ee
where $w:=(v_m(c)-v)^{1/2}$ (cf.\ \cite{H}).  The first derivative of $m$ is 
\be\nn\begin{aligned}
m'(c) & = 4\int \frac{d}{dc}\Bigl\{ \left(\kappa(v_m(c)-w^2)\right)^{1/2}
\left(-2F(v_m(c)-w^2,c) \right)^{1/2}\Bigr\}\,w\;dw\\
&= 4\int
\frac{\kappa_v(v_m(c)-w^2)v_m'(c)}{2\left(\kappa(v_m(c)-w^2)\right)^{1/2}}
\left(-2F(v_m(c)-w^2,c) \right)^{1/2}\,w\;dw\;\; + \\
&\quad + 4 \int \frac{\left(\kappa(v_m(c)-w^2)\right)^{1/2} \left(-
F_v(v_m(c)-w^2,c)v_m'(c) - F_c(v_m(c)-w^2,c)\right)  }{\left(-2F(v_m(c)-w^2,c)
\right)^{1/2}}\;w\;dw     
 \end{aligned}\ee
which, due to
\be\nn
\frac{\partial}{\partial v} \left(\left(\kappa(v) (-2F(v,c))\right)^{1/2}\right)
= \frac{\kappa_v(v)(-2F(v,c))^{1/2}}{2(\kappa(v))^{1/2}} -
\frac{\kappa(v)^{1/2}F_v(v,c)}{(-2F(v,c))^{1/2}},
\ee
simplifies to 
\be\nn
m'(c) = -4 \int
\frac{\left(\kappa(v_m(c)-w^2)\right)^{1/2}F_c(v_m(c)-w^2,c)}{
\left(-2F(v_m(c)-w^2,c) \right)^{1/2}}\;w\;dw.
\ee
The second derivative of $m$ can be written in the form  
\be\nn
\begin{aligned}
 m''(c)&=2 \int_{v_*}^{v_m(c)} \frac{A(v,c) +
B(v,c)}{\left(\kappa(v)\right)^{1/2}\left(-2F(v,c)\right)^{3/2}}\;dv
\end{aligned}
\ee
with 
$A(v,c)= F(v,c)\kappa_v(v)v_m'(c)F_c(v,c)
$
and
$$
B(v,c) = \kappa(v)\left( (v-v_*)\left( 2F(v,c)\left((v-v_*)+2cv_m'(c)\right)
-c(v-v_*)\left( F_v(v,c)v_m'(c)+F_c(v,c) \right) \right)  \right).
$$
As $F_c(v,0)=0=A(v,0)$ 
and
\be\nn
B(v,0)= 2\kappa(v)(v-v_*)^2F(v,0)<0\quad \text{for all }v\in (v_*,v_m(0)),
\ee
we indeed have $m''(0)<0.$
\\

{\bf Remark.} Certain of the results by Zumbrun \cite{Z} on 
the Bona-Sachs model ($\kappa(V)\equiv 1$ and $p(V)=-V+V^q$ with $q\geq 2$) 
and by De Bouard \cite{DeB} on a Gross-Pitaevskii model ($\kappa(V) =1/(4V^{4})$ and $p(V)=\alpha/V^2 - \beta/V^3$ with $\alpha,\beta>0$)
appear as interesting special cases of Theorem 1.
\\

I obtained parts of this result during my doctoral dissertation under the supervision of Heinrich Freist\"uhler; I am very grateful to him for helpful discussions.

\end{document}